\documentclass[12pt,reqno]{amsart}

\usepackage{amsmath, epsfig, cite,color}
\usepackage{amssymb}
\usepackage{amsfonts}
\usepackage{latexsym}
\usepackage{graphicx}
\usepackage{algorithmic,algorithm}
\usepackage{amsthm,array,multirow,longtable}
\usepackage{tikz,caption}
\usetikzlibrary{decorations.pathreplacing}

\newtheorem{theorem}{Theorem}

\numberwithin{equation}{section}

\newcommand{\nc}{\newcommand}
\nc{\sha}{\mbox{\cyr X}}  
\font\cyr=wncyr10

\makeatletter \@addtoreset{equation}{section} \makeatother

\begin{document}
\allowdisplaybreaks
\title[One level summations for powers of Fibonacci polynomials]
{One level summations for powers of Fibonacci and Lucas polynomials}
\author[H. Prodinger]{Helmut Prodinger}
\address[H. Prodinger]{Department of Mathematics, University of Stellenbosch,
7602 Stellenbosch, South Africa} \email{hproding@sun.ac.za}



\date{\today}

\begin{abstract} Powers of Fibonacci polynomials are expressed as single sums, improving on a double sum recently seen in the literature.  

\end{abstract}

\maketitle

\section{Introduction}

The generating function of Fibonacci numbers is
\begin{equation*}
\sum_{n\ge0}F_nz^n=\frac{z}{1-z-z^2}=z\sum_{k\ge0}(z+z^2)^k=\sum_{k\ge0}z^{k+1}(1+z)^k.
\end{equation*}
Comparing coefficients on both sides,
\begin{equation*}
F_{n+1}=\sum_{k\ge0}[z^{n-k}](1+z)^k=\sum_{0\le k\le n}\binom{k}{n-k}=\sum_{0\le i\le n/2}\binom{n-i}{i}.
\end{equation*}
This is of course a very classical formula \cite{wiki}.

Recently, the following formula was found \cite{Thai2}:
\begin{equation*}
F_{n+1}^2=\sum_{0\le j \le i\le 2n/3}\binom{i}{j}\binom{2n-2i-j}{i},
\end{equation*}
which contains already \emph{two} levels of summation.

In \cite{Fron}, the author derived, among many other things, the formula
\begin{align*}
	F_n^2&=\frac25F_{2(n+1)}-\frac35F_{2n}-\frac25(-1)^n=\frac{L_{2n}}{5}+\frac{2(-1)^{n-1}}{5},
\end{align*}
with Lucas numbers, satisfying
\begin{equation*}
	L_n=\sum_{0 \le k \le n/2}\frac{n}{n-k}\binom{n-k}{k},\quad n\ge1.
\end{equation*}

Consequently, we find for $n\ge1$
\begin{align*}
	F_n^2&=\frac15 \sum_{0 \le k \le n}\frac{2n}{2n-k}\binom{2n-k}{k}+\frac{2(-1)^{n-1}}{5},
\end{align*}
this result has only \emph{one} level of summations. For fun, we go one step further. Since \cite{Fron}
\begin{equation*}
F_n^3=\frac15F_{3n}-\frac35(-1)^nF_n,
\end{equation*}
we find
\begin{equation*}
	F_n^3=\frac15 \sum_{0\le i\le (3n-1)/2}\binom{3n-1-i}{i} -\frac35(-1)^n\sum_{0\le i\le (n-1)/2}\binom{n-1-i}{i},
\end{equation*}
again only one level of summation.

The rest of the paper is devoted to  Fibonacci polynomials, defined by
\begin{equation*}
F_{n+2}=xF_{n+1}+F_{n},\quad F_0=0,\ F_1=1,
\end{equation*}
and \emph{arbitrary} powers $F_n^d$, for positive integers $d$, not just squares.
There are also Lucas polynomials
\begin{equation*}
	L_{n+2}=xL_{n+1}+L_{n},\quad L_0=2,\ L_1=x,
\end{equation*}
and if one understands these two, one understand arbitrary initial conditions, as well a the recursion $f_{n+2}=xf_{n+1}+yf_{n}$, since one parameter can be
eliminated by a simple substitution.

From now on $F$ and $L$ mean Fibonacci polynomials (or generalized Fibonacci numbers) and, respectively, Lucas polynomials (or generalized Lucas numbers).

\section{Linearizing powers of Fibonacci polynomials and Lucas polynomials}

First, solving $\lambda^2=x\lambda+1$, with the two roots
\begin{equation*}
	\alpha,\beta=\frac{x\pm\sqrt{x^2+4}}{2},
\end{equation*}
  we have the \emph{Binet forms}
\begin{equation*}
F_n=\frac{\alpha^n-\beta^n}{\alpha-\beta} \quad\text{and}\quad L_n=\alpha^n+\beta^n.
\end{equation*}
Note that $\alpha-\beta=\sqrt{x^2+1}$ and $\alpha \beta=-1$. 
The expansions for $n\ge1$ are \cite{wiki}
\begin{align*}
F_n&=\sum_{0\le k\le (n-1)/2}\binom{n-1-k}{k}x^{n-1-2k},\\
L_n&=\sum_{0\le k\le n/2}\frac{n}{n-k}\binom{n-k}{k}x^{n-2k}.
\end{align*}
It is to be noted that these expansion only contain one level of summation.

Doing a few experiments, we found  two sets of formulae.
\begin{theorem}
	
For $n\ge0$ and $d\ge0$,
\begin{equation*}
F_n^{2d+1}=\frac1{(x^2+4)^d}\sum_{0\le s \le d}(-1)^{s(n-1)}\binom{2d+1}{s}F_{(2d+1-2s)n}.
\end{equation*}
For $n\ge0$ and $d\ge1$,
\begin{equation*}
	F_n^{2d}=\frac1{(x^2+4)^d}\sum_{0\le s \le d-1}(-1)^{s(n-1)}\binom{2d}{s}L_{(2d-2s)n}+\frac1{(x^2+4)^d}\binom{2d}{d}(-1)^{d(n-1)}.
\end{equation*}
\end{theorem}

We will provide two sample computations that illustrate the general situation:
\begin{align*}
F_n^7&=\frac{(\alpha^n-\beta^n)^7}{(\alpha-\beta)^7}\\*
&=\frac{\alpha^{7n}-7\alpha^{6n}\beta^{n}+21\alpha^{5n}\beta^{2n}-35\alpha^{4n}\beta^{3n} +35\alpha^{3n}\beta^{4n}
-21\alpha^{2n}\beta^{5n}+7\alpha^{n}\beta^{6n}-\beta^{7n} }{(\alpha-\beta)(x^2+4)^3}\\
&=\frac{\alpha^{7n}-\beta^{7n}-7(-1)^n\alpha^{5n}+7(-1)^n\beta^{5n}+21\alpha^{3n}	-21\beta^{3n} -35(-1)^n\alpha +35(-1)^n\beta
}{(\alpha-\beta)(x^2+4)^3}\\
&=\frac{F_{7n}-7(-1)^nF_{5n}+21F_{3n} -35(-1)^nF_n}{(x^2+4)^3}.
\end{align*}
The situation for even exponents is similar, but there is a middle term with no matched term:
\begin{align*}
	F_n^6&=\frac{(\alpha^n-\beta^n)^6}{(\alpha-\beta)^6}\\*
	&=\frac{\alpha^{6n}-6\alpha^{5n}\beta^{n}+15\alpha^{4n}\beta^{2n}-20\alpha^{3n}\beta^{3n} +15\alpha^{2n}\beta^{4n}
		-6\alpha^{n}\beta^{5n}+\beta^{6n} }{(x^2+4)^3}\\
	&=\frac{\alpha^{6n}-6(-1)^n\alpha^{4n}+15\alpha^{2n}-20(-1)^n+15\beta^{2n}
		-6(-1)^n\beta^{4n}+\beta^{6n} }{(x^2+4)^3}\\
	&=\frac{L_{6n}-6(-1)^nL_{4n}+15L_{2n}-20(-1)^n		}{(x^2+4)^3}.
\end{align*}
\textsc{Remark.} In \cite{Fron}, a similar approach was used for Fibonacci resp.\ Lucas \emph{numbers} (not polynomials). Since the notion of Lucas numbers was not
used in the even instance, the description was necessarily clumsier than here.

Now we move to powers of Lucas polynomials, which are a bit simpler to describe, although it is still better to distinguish odd resp.\ even exponents:
\begin{align*}
L_n^{2d+1}&=\sum_{0\le s \le 2d+1}\binom{2d+1}s\alpha^{(2d+1-s)n}\beta^{sn}\\
&=\sum_{0\le s \le d}\binom{2d+1}s(-1)^{sn}\Big(\alpha^{(2d+1-2s)n}+\beta^{(2d+1-2s)n}\Big)
\\&=\sum_{0\le s \le d}\binom{2d+1}s(-1)^{sn}L_{(2d+1-2s)n}.
\end{align*}
\begin{align*}
	L_n^{2d}&=\sum_{0\le s \le 2d}\binom{2d}s\alpha^{(2d-s)n}\beta^{sn}\\*
		&=\sum_{0\le s \le d-1}\binom{2d}s(-1)^{sn}\Big(\alpha^{(2d-2s)n}+\beta^{(2d-2s)n}\Big)+\binom{2d}{d}(-1)^{dn}
	\\&=\sum_{0\le s \le d-1}\binom{2d}s(-1)^{sn}L_{(2d+1-2s)n}+\binom{2d}{d}(-1)^{dn}.
\end{align*}

\section{Conclusion}

The authors in \cite{Thai1, Thai2} were also interested in writing $F_nF_{n+1}$ as a sum. We can compute this more generally:
\begin{equation*} 
F_nF_{n+d}=\frac1{x^2+4}L_{2n+d}-\frac{L_d}{x^2+4}(-1)^n,
\end{equation*}
 and this is clearly a single sum (just one level of summations). The justification of the formula goes via the Binet forms of the product and expanding.

\end{document}